\documentclass[12pt]{article}
\usepackage{amsfonts}
\usepackage{amsfonts,mathrsfs}
\usepackage{CJK,fancyhdr,amscd}
\usepackage{anysize}
\usepackage{graphicx,epsfig}
\usepackage{amsthm,latexsym,amsmath,amssymb}
\usepackage[perpage,symbol*]{footmisc}

\newtheorem{thm}{Theorem}[section]
\newtheorem{defi}[thm]{Definition}
\newtheorem{prop}[thm]{Proposition}
\newtheorem{lem}[thm]{Lemma}
\newtheorem{cor}[thm]{Corollary}
\newtheorem{eg}{Example}[section]
\newtheorem{rem}[thm]{Remark}

\makeatletter \@addtoreset{equation}{section} \makeatother
\begin{document}
\baselineskip=20pt  \hoffset=-3cm \voffset=0cm \oddsidemargin=3.2cm
\evensidemargin=3.2cm \thispagestyle{empty}\vspace{10cm}
\title{\textbf{Symmetrical Symplectic Capacity with Applications}}
\author{   Chungen
Liu\thanks{Corresponding author: liucg@nankai.edu.cn.}
\thanks{ Partially  supported by NNSF of
China(10531050,10621101) and 973 Program of
STM(2006CB805903).}\;\;\;\;\;\;\;\;Qi Wang\thanks{Partially
supported by the NNSF of
China(10701043).eliada@mail.nankai.edu.cn.}\\
    \\School of Mathematics and LPMC \\ Nankai
University, Tianjin 300071, P.R.China}
\date{} \maketitle

 \noindent{\bf Abstract:} {\small  In this paper, we first introduce the concept of symmetrical
symplectic capacity for  symmetrical symplectic manifolds, and by
using this symmetrical symplectic capacity theory  we prove that
there exists at least one  symmetric closed characteristic (brake
orbit and $S$-invariant brake orbit are two examples) on prescribed
symmetric energy surface which has a compact neighborhood with
finite symmetrical symplectic capacity.}



\section{Introduction and main results}
The purpose of this paper is to study the existence of the symmetric
periodic solutions of Hamiltonian systems in the presence of
symmetry for the manifold and also for the Hamiltonian functions. A
very famous example is the figure-eight orbit in planar three-body
problem with equal masses(see \cite{CM}). It is the orbit with two
different symmetries: cyclic symmetry and generalized brake
symmetry. In this paper we consider the existence of symmetric
orbits of smooth Hamiltonian systems with some symmetries. An
important case is the existence of brake orbits  on the manifolds
with the brake symmetry. For this purpose, we first study the
symmetrical symplectic capacity theory for the symplectic manifolds
with corresponding symmetry.

 Symplectic capacity is an  important symplectic
invariant. It was first discovered by I.Ekeland and H.Hofer in
\cite{E-H} and \cite{E-H2} for subsets of $\mathbb{R}^{2n}$ in their
search for periodic solutions of Hamiltonian systems on fixed energy
surfaces. We call it the Ekeland-Hofer capacity and denote it by
$c_{EH}$. This concept was extended to general symplectic manifolds
by H.Hofer and E.Zehnder in \cite{H-Z2} and \cite{H-Z}. We call it
the Hofer-Zehnder capacity and denote it by $c_0$. As  examples of
symplectic capacity, Gromov's width $\mathcal{W}_G$ defined in
(\cite{Gr}) is the smallest symplectic capacity, Hofer's
displacement energy $d$ defined in (\cite{H}) is also a symplectic
capacity, the Floer-Hofer capacity $c_{FH}$ defined in (\cite{He})
which can be viewed as a variant of Ekeland-Hofer capacity $c_{EH}$,
and Viterbo's generating function capacity $c_V$ defined in
(\cite{V}) is also a symplectic capacity. The symplectic capacities
were applied to the study of many symplectic topology problems, see
\cite{H-Z}, \cite{McSa}, \cite{V2} and the references therein for
more details.

In this paper, we introduce a symmetrical capacity on some
symmetrical symplectic  manifolds.  For the brake symmetry case,
we say that a symplectic manifold $(\mathcal{M},\omega)$ is brake
 symmetrical ($\varphi$-symmetric) if there is an antisymplectic involution $\varphi:\mathcal{M}\to
 \mathcal{M}$ satisfying $\varphi^2={\rm id}$, $\varphi^*\omega=-\omega$ and the fixed point set ${\rm Fix}(\varphi)\neq
 \emptyset$. It is well
 known that the fixed point set $\mathcal{L}$ of $\varphi$ is a Lagrangian
 submaifold of $\mathcal{M}$ if it is not empty.  We denote the
 $\varphi$-symmetric symplectic manifold $\mathcal{M}$ by  $(\mathcal{M},\mathcal{L},\omega,\varphi
 )$. For a $\varphi$-symmetric symplectic manifold $(\mathcal{M},\mathcal{L},\omega,\varphi
 )$, in this paper we first develop a
 {symmetrical symplectic capacity} $c_\varphi(\mathcal{M})$ in subsection 2.1. When a symplectic manifold
 $(M,\omega)$ is provided with two different symmetries, for
 example,
 the $\varphi$-symmetry and
  a cyclic symmetry $S$, for some special cases we also
 introduce a capacity $c_{\varphi,S}(\mathcal{M})$.  For example, in $(\mathbb{R}^{2n},\omega_0)$, we  choose
 $\varphi$ as a linear mapping $N_0:\mathbb{R}^{2n}\to\mathbb{R}^{2n}$ with
 $N_0=\left(\begin{matrix}-I_n & 0\\0 & I_n\end{matrix}\right)$, and
 an orthogonal symplectic matrix $S$
 with $S^m={\rm id}$ for some $2\le m\in \mathbb{N}$, we introduce
 a symmetrical
 capacity $c_{N_0,S}(U)$ for $(N_0,S)$ invariant subset $U$ of $\mathbb{R}^{2n}$ in subsection 3.1.

We note that for a general symplectic manifold it is not easy to
determine the finiteness of its Hofer-Zehnder's capacity. There are
few  results about the finiteness of
 symplectic capacity for some
 special symplectic manifolds(see for example
\cite{H-V},\cite{Jiang},\cite{Lu},\cite{M}). It's also difficult
for us to prove the finiteness of symmetrical symplectic capacity
in general, in Theorem \ref{main-result} and Lemma \ref{Tn} below
we give some special examples with finite symmetrical symplectic
capacity.

As the applications of the symmetrical symplectic capacity
$c_\varphi$,  we consider the existence of brake orbits (see
Definition \ref{2.2} below) on  energy hypersurfaces in symmetrical
symplectic manifolds. The main results read as follows.

\begin{thm}\label{H-Z} For a $\varphi$-symmetric symplectic manifold
$(\mathcal{M},\mathcal{L},\omega,\varphi
 )$,
let $\Sigma=H^{-1}(1)$ be the compact  regular energy surface of a
$\varphi$-invariant Hamiltonian function $H\in
C^2(\mathcal{M},\mathbb{R})$. Suppose $\Sigma\cap
\mathcal{L}\neq\emptyset$ and there is an open neighborhood $U$ of
$\Sigma$ such that $c_{\varphi}(U,\omega)<\infty$. Then there exists
a sequence $\lambda_j\rightarrow 1$, $j\to +\infty$, such that every
energy surface $\Sigma_{\lambda_j}=H^{-1}(\lambda_j)$ possesses a
brake orbit of the Hamiltonian vector field $X_H$.
\end{thm}

For a $\varphi$-symmetric contact-type hypersurface $\Sigma\in
\mathcal{S}_{\varphi}$ in $(\mathcal{M},\mathcal{L},\omega,\varphi)$
which are defined in Definition \ref{3.4} and \ref{r-c-t} below, we
have
\begin{thm}\label{C.Viterbo}
Suppose that the $\varphi$-contact type hypersurface $\Sigma\in
\mathcal{S}_{\varphi}$ has a $\varphi$-invariant neighborhood $U$
with $c_\varphi(U,\omega)<\infty$, then $\Sigma$ possesses a closed
brake-characteristic.
\end{thm}

We note that for a compact $\varphi$-contact type hypersurface
$\Sigma$ in $(\mathbb{R}^{2n},\omega_0)$ with $\varphi=N_0$, it is
clear that $\Sigma$ has a $\varphi$-invariant neighborhood $U$ with
$c_\varphi(U,\omega)<\infty$. So $\Sigma$ always possesses a closed
brake-characteristic.

In section 3.2,  as  applications of $c_{N_0,S}$, we consider the
existence of $S$-symmetrical brake orbits on energy hypersurface in
$(\mathbb{R}^{2n}, \omega_0)$, and get the following result.
\begin{thm}\label{T1.3}
Let $\Sigma=H^{-1}(1)$ be the compact  regular energy surface of a
$(N_0, S)$-invariant function $H\in
C^2(\mathbb{R}^{2n},\mathbb{R})$. Suppose $\Sigma$ is the boundary
of a bounded domain $O$ in $\mathbb{R}^{2n}$ with $0\in O$. Then
there is an open neighborhood $U$ of $\Sigma$ such that
$c_{N_0,S}(U)<\infty$ and there exists a sequence
$\lambda_j\rightarrow 1$,$j\to +\infty$, such that every energy
surface $\Sigma_{\lambda_j}=H^{-1}(\lambda_j)$ possesses a
$S$-symmetrical brake orbit of the Hamiltonian vector field $X_H$.
\end{thm}

 We shall note that on a fixed energy surface there may be no closed
 characteristic(see \cite{M.Her1}, \cite{M.Her2} for counter
 examples). But for
the case of $(\mathbb{R}^{2n},\omega_0)$ as considered in Example
\ref{eg1.1} below, if the the $N_0$-invariant hypersurface
$\Sigma=H^{-1}(1)$ is star-shaped,  Rabinowitz in 1987 \cite{R}
proved that if $x\cdot H'(x)\neq 0$ for all $x\in \Sigma$ then there
exist at least one brake orbits on $\Sigma$, which has been
generalized by Corollary \ref{3.7} below in this paper. If the
$N_0$-invariant hypersurface $H^{-1}(1)$ is $\sqrt2$-pinched, A.
Szulkin in 1989 \cite{Sz} proved that it possesses at least $n$
geometrically distinct brake orbits. If the $N_0$-invariant
hypersurface $\Sigma=H^{-1}(1)$ is convex and central symmetric,
that is $\Sigma=-\Sigma$, Y. Long, D. Zhang and C. Zhu in 2006
\cite{Y-D-Z} proved that $\Sigma$ possesses at least two
geometrically distinct brake orbits. Recently, in 2009 \cite{L-Z},
D. Zhang and the first author of this paper proved that a convex and
central symmetric hypersurface $\Sigma\subset \mathbb{R}^{2n}$
possesses at least $[\frac{n}{2}]+1$ geometrically distinct brake
orbits, and if all brake orbits on $\Sigma$ are nondegenerate, then
$\Sigma$ possesses at least $n$ geometrically distinct brake orbits.

 For brake boundary value problems of non-autonomous Hamiltonian, one can
 refer the papers
\cite {L-L}, \cite{W-L-L} and \cite{Z}.  For the existence and
multiplicity of closed characteristics on prescribed energy surface,
one can further refer the papers \cite{LLZ,LZ,R2,R3,St,V3,Wei1,Wei2}
and the references therein.

\section{Symmetrical Symplectic Capacity and Its Applications }
In this section, we first introduce the concept of symmetric
symplectic capacity and develop some properties for this kind of
capacity. As  applications, we then prove Theorem \ref{H-Z} and Theorem
\ref{C.Viterbo}.
\subsection{Symmetrical Symplectic Capacity}\label{2.1}
 \begin{defi}\label{def1.1} A  symplectic manifold $(\mathcal{M},\omega)$ is called a symmetrical
symplectic manifold, if there exists a diffeomorphism $\varphi:
\mathcal{M}\rightarrow \mathcal{M}$, and a Lagrangian submanifold
$\mathcal{L}$  of $\mathcal{M}$ satisfying
\begin{equation}\label{1.1}\mathcal{L}=Fix(\varphi),\; \varphi^2=id|_{\mathcal{M}}\;{\rm
and}\;\varphi^*\omega=-\omega.\end{equation}
\end{defi}
 From now on, we always denote by $(\mathcal{M},\mathcal{L},\omega,\varphi)$
the { symmetrical symplectic manifold} with $\varphi$ and
$\mathcal{L}$ satisfying condition \eqref{1.1}. For symmetrical
symplectic manifolds,  we have the following examples.
\begin{eg}\label{eg1.1}
 {\rm The linear symplectic space $(\mathbb{R}^{2n},\omega_0)$ with $\omega_0=\displaystyle
 \sum_{k=1}^ndx_k\wedge dy_k$, let $N\in\mathcal
{L}(\mathbb{R}^{2n})$ satisfying the following conditions
\begin{equation}\label{N}
 N^TJN=-J,\;\;\;\; N^2=I_{2n\times 2n},
\end{equation}
 where $J=\left(\begin{matrix}0 & -I_{n\times n}\\I_{n\times n} &
 0\end{matrix}\right)$. It is easy to see that $L_N:=ker(N-I_{2n\times 2n})$
 is an $n$ dimensional Lagrangian subspace of $\mathbb{R}^{2n}$,
 and $(\mathbb{R}^{2n}, L_N, \omega_0, N)$ satisfies the conditions of Definition
 \ref{def1.1},
 in particular, for $N=N_0=\left(\begin{matrix}-I_{n\times n} &0\\0
 &I_{n\times n}\end{matrix}\right)$, in this case $L_N=L_0=\{0\}\times
 \mathbb{R}^{n}$.
 }
\end{eg}
\begin{eg}\label{eg1.2}
{\rm
 Let $\mathcal{N}$ be an $n$ dimensional smooth manifold, and $\mathcal{M}=T^*\mathcal{N}$.  $(\mathcal{M},
 \omega)$ is a symplectic manifold with $\omega$ being the canonical
 symplectic form.
 Let $\varphi: \mathcal{M}\rightarrow\mathcal{M}$,
 $\varphi(x,\xi)=(x,-\xi),\; \forall \,(x,\xi)\in\mathcal{M}$, then ${\rm Fix}(\varphi)=\mathcal{N}$ and $(\mathcal{M},\mathcal{N}, \omega,
 \varphi)$ is a symmetrical symplectic manifold.
 }
\end{eg}
\begin{eg}\label{eg1.3}
{\rm In the case of $\mathcal{N}=T^n:=\mathbb{R}^n/2\pi \mathbb{Z}$
as in Example \ref{eg1.2}, but the involution $\varphi$ is defined
by $\varphi:T^*(T^n)\rightarrow T^*(T^n)$,
$\varphi(\theta_1,\cdots,\theta_n,\xi_1,\cdots,\xi_n)=(-\theta_1,\cdots,-\theta_n,\xi_1,\cdots,\xi_n)$,
so ${\rm Fix}(\varphi)=\mathcal{L}:=(0,\cdots,0)\times
\mathbb{R}^n\cup(\pi,\cdots,\pi)\times\mathbb{R}^n$, and
$(T^*(T^n),\mathcal{L},\omega,\varphi)$ is a symmetrical symplectic
manifold.

}
\end{eg}
 \noindent For a given symmetrical symplectic manifold $(\mathcal{M},\mathcal{L},\omega,\varphi)$,
  we denote by $\mathcal{H}(\mathcal{M},\mathcal{L},\omega,\varphi)$ the
set of $C^2$  functions $H:\mathcal{M}\to \mathbb{R}$ satisfying the following four properties:\\
 (H1)  There is a compact set $K\subset \mathcal{M}$ (depending on $H$ ) such that
 $K\subset(\mathcal{M}\setminus\partial \mathcal{M})$ and
 $$H(\mathcal{M}\setminus K)\equiv m(H) \;\;\;\text{(a constant)}.$$
 (H2)  There is an open set $\mathcal {O}\subset \mathcal{M}\;(\text{depending on }H) \;{\rm and }\;\mathcal{O}\cap\mathcal{L}\neq\emptyset\;$
  on which
 $$H(\mathcal {O})\equiv 0.$$
 (H3)   $0\leq H(x)\leq m(H)$ for all $x\in \mathcal{M}$.\\
 (H4)  $H({\varphi}(x))=H(x).$

\begin{defi}\label{2.2}Suppose $H\in  C^2(\mathcal{M}, \mathbb{R})$
satisfying {\rm (H4)},
 consider the
Hamiltonian system on $\mathcal{M}$
\begin{equation}\label{brake-orbit}
 \left\{ \begin {array}{lll}
\dot{x}(t)=X_H(x(t)),\\
 x(-t)=\varphi(x(t)),\\
x(T+t)=x(t),
  \end{array}\right.
 \end{equation}
where $X_H$ is the Hamiltonian vector field with respect to the
function $H$.
  A solution $(x,T)$ of \eqref{brake-orbit} is called a  brake
orbit of  the Hamiltonian vector field $X_H$.
\end{defi}
 A function $H\in\mathcal{H}(\mathcal{M},\mathcal{L},\omega,\varphi)$ is called admissible if \eqref{brake-orbit} has no brake orbits, or all
the brake orbits  of \eqref{brake-orbit} on $\mathcal {M}$
are either constant, ie., $x(t)\equiv x(0)$, $\forall t\in\mathbb{R}$ or have the minimal
period $T>1$. Denote the set of admissible functions by
$\mathcal{H}_a(\mathcal{M},\mathcal{L},\omega,\varphi)\subseteq\mathcal{H}(\mathcal{M},\mathcal{L},\omega,\varphi)$.
\begin{defi}\label{cc} We define the symmetrical symplectic capacity on a
symmetrical symplectic manifold
$(\mathcal{M},\mathcal{L},\omega,\varphi)$ by
\begin{equation}\label{c}
c_{\varphi}(\mathcal{M},\omega)=\sup\{m(H)|H\in\mathcal
{H}_a(\mathcal{M},\mathcal{L},\omega,\varphi)\}.
\end{equation}
\end{defi}
 \noindent From the definition, we see that
 if $c_{\varphi}(\mathcal{M},\omega)<\infty$,
for every function $H$ in
$\mathcal{H}(\mathcal{M},\mathcal{L},\omega,\varphi)$  satisfying
$m(H)>c_{\varphi}(\mathcal{M},\omega)$, the vector field $X_H$
possesses a nonconstant  brake orbit with minimal period $0<T\leq
1$, and $c_{\varphi}(\mathcal{M},\omega)$ is the infimum
of the real numbers having this property.\\

 \begin{thm}\label{main-result}
 $c_{\varphi}$ has the following properties: \\
(A) $c_{\varphi_1}(\mathcal{M}_1,\omega_1)\leq
c_{\varphi_2}(\mathcal{M}_2,\omega_2)$ provided there exists a
symplectic embedding $\phi:\mathcal{M}_1\rightarrow \mathcal{M}_2$
satisfying $\phi(\varphi_1(z))=\varphi_2(\phi(z))$ for $\forall z\in
\mathcal{M}_1$, here $(\mathcal{M}_1,\omega_1)$ and
$(\mathcal{M}_2,\omega_2)$ are two
symmetrical symplectic manifolds with same dimension.\\
(B) $c_{\varphi}(
\mathcal{M},\alpha\omega)=|\alpha|c_{\varphi}(\mathcal{M},\omega)$,
$\alpha\neq 0$, $\alpha\in \mathbb{R}$.\\
(C) $c_{N}(B(1),\omega_0)=\pi$, where $B(1)$ is the unit ball in
$\mathbb{R}^{2n}$,  $N$ and $\omega_0$ are defined in Example
\ref{1.1} and $N=N^T$.
\end{thm}
\noindent From Theorem \ref{main-result}, we see that the symmetric
symplectic capacity $c_{\varphi}(\mathcal{M},\omega)$ satisfies all
properties of the general symplectic capacity in the sense of
symmetric category(c.f., \cite{E-H,E-H2,H-Z2, H-Z}). The proof of
Theorem \ref{main-result} is similar to that as in
\cite{E-H,E-H2,H-Z2, H-Z}. We complete the proof of Theorem
\ref{main-result} via the following lemmas.

\begin{lem}\label{A}
$c_{\varphi}$ satisfies the  properties (A) and (B).
\end{lem}
\noindent {\it Proof.} We divided the proof into two steps.\\ {\it
Step 1. Proof of property (A).}  Define a map $\phi_*:\mathcal
{H}(\mathcal{M}_1,\mathcal{L}_1,\omega_1,\varphi_1)\rightarrow\mathcal
{H}(\mathcal{M}_2,\mathcal{L}_2,\omega_2,\varphi_2)$ by
$$
\phi_*(H)=\left\{\begin{array}{ll} H\circ\phi^{-1}(x)
\;\;&\text{if}\;\;x\in \phi(\mathcal{M}_1),\\
m(H)\;\;&\text{if}\;\;x\notin \phi(\mathcal{M}_1).
\end{array}\right.
$$
Note that if $K\subset \mathcal{M}_1\setminus  \partial
\mathcal{M}_1$ for a compact set $K\subset \mathcal{M}_1$, then also
$\phi(K)\subset \mathcal{M}_2\setminus
\partial \mathcal{M}_2$, so there holds $\phi_*(H)\in \mathcal
{H}(\mathcal{M}_2,\mathcal{L}_2,\omega_2,\varphi_2)$. Clearly
$m(\phi_*(H))=m(H)$. Since $\phi$ is symplectic and satisfies
$\phi(\varphi(z))=\varphi(\phi(z))$ for $z\in \mathcal{M}_1$, we
have $\phi_*(\mathcal
{H}_a(\mathcal{M}_1,\mathcal{L}_1,\omega_1,\varphi_1))\subset
\mathcal {H}_a(\mathcal{M}_2,\mathcal{L}_2,\omega_2,\varphi_2)$. This implies the property (A).\\
{\it Step 2. Proof of property (B).} Assume $\alpha\neq 0$ and
define the bijection $\psi$:
$\mathcal{H}(\mathcal{M},\mathcal{L},\omega,\varphi)\rightarrow\mathcal{H}(\mathcal{M},\mathcal{L},\omega,\varphi)$
by $\psi:H\mapsto H_\alpha:=|\alpha|H$. Clearly
$m(H_\alpha)=|\alpha|m(H)$, By the definition of $X_H$, we have
$\frac{\alpha}{|\alpha|}X_{H_\alpha}=X_{H}$ on $\mathcal{M}$.
Therefor, $X_{H_\alpha}$ and $X_H$ have the same brake orbits with
the same periods. It implies that $\psi$ is also a bijection between
$\mathcal {H}_a(\mathcal{M},\mathcal{L},\omega,\varphi)$ and
$\mathcal {H}_a(\mathcal{M},\mathcal{L},\alpha\omega,\varphi)$. Thus
the property (B) is true. $\hfill\Box$

 For the proof of
property (C), we note that it is enough to prove it for $N=N_0$. In
fact,  there exists an orthogonal symplectic matrix $P$ satisfying
$$P^{-1}NP=N_0.$$
That is to say $P:(B(1),L_N\cap B(1),\omega,N)\rightarrow
(B(1),L_0\cap B(1),\omega,N_0)$ is a symplectic diffeomorphism
satisfying the condition in property (A), so  we have
$$c_N(B(1))=c_{N_0}(B(1)).$$
A brake orbit is naturally a periodic orbit, from the definitions of
$c_0$ and $c_\varphi$, there holds
$c_\varphi(\mathcal{M},\omega)\geq c_0(\mathcal{M},\omega)$. in
particular, we have
$$c_{N_0}(B(1))\geq c_0(B(1))= \pi.$$
In order to prove property (C), we need to  prove
$c_{N_0}(B(1))\leq\pi$. By property (A),  it is enough to  prove
\begin{equation}\label{z(1)}c_{N_0}(Z(1), \omega_0)\leq\pi, \end{equation} where
$Z(1)=\left\{(x,y)\in\mathbb{R}^{2n}|x_1^2+y_1^2<1\right\}$.

From now on we  assume that the Hamiltonian function $H\in \mathcal
{H}(Z(1),L_0\cap Z(1),\omega_0,N_0)$ (denote by $\mathcal{H}(Z(1))$
for short in the sequel) satisfying $m(H)>\pi$. The remaining  of
this subsection is to prove that the vector field $X_H$ possesses a
nonconstant  brake orbit with minimal period $0<T\leq 1$. First, we
have

\begin{lem}\label{L2.6} Suppose the Hamiltonian function $H\in \mathcal
{H}(Z(1))$ satisfying $m(H)>\pi$. Then there is a compactly
supported symplectic diffeomorphism of $Z(1)$, $\psi:Z(1)\to Z(1)$,
i.e., the closure of the set $\{x|\psi(x)\neq x\}$ is a compact
subset of $Z(1)$, satisfying  $\psi(N_0 z)=N_0\psi(z),\;\forall z\in
Z(1)$, and $H\circ\psi:Z(1)\to \mathbb{R}$ vanishes in an open
neighborhood of the origin $0$.
\end{lem}
\noindent{\it Proof.} Since $H\in\mathcal{H}(Z(1))$, there is an
open set $\mathcal {O}\subset Z(1)\;{\rm and }\;\mathcal{O}\cap
L_0\neq\emptyset\;$
  on which
 $H(\mathcal {O})\equiv 0.$
  If $\{0\}\in\mathcal{O}$, since $H(\mathcal{O})\equiv 0$, $H$
vanishes in an open neighborhood of the origin.  In this case, we
define $\psi={\rm id}|_{Z(1)}$.
 Otherwise since
$\mathcal{O}\cap L_0\neq\emptyset$, by choosing a point $z_0\in
\mathcal{O}\cap L_0$, $z_0\neq 0$, we have the orthogonal
decomposition $\mathbb{R}^{2n}=span\{z_0\}\oplus\mathbb{R}^{2n-1}$,
$z=(tz_0,z_{n-1})\in \mathbb{R}^{2n}$. Denote by
$A(\delta)=(-\delta,1+\delta)\times B_{\mathbb{R}^{2n-1}}(\delta)$,
with $\delta>0$  small enough such that $A(2\delta)\subset Z(1)$.
Denote by $\chi_{A(\delta)}(z)$ the characteristic function of
$A(\delta)$. That is to say
\begin{equation}
\chi_{A(\delta)}(z)=\left\{\begin{array}{ll} 1,\;\;&{\rm if}\; z\in A(\delta),\\
0, &{\rm otherwise}.\end{array} \right.
\end{equation}
Choosing  a smooth function $\sigma_\delta\in
C^\infty_0(\mathbb{R}^{2n}, \mathbb{R}^+)$  with compact support in
$B(\delta/2)\subset\mathbb{R}^{2n}$ and
$\int_{\mathbb{R}^{2n}}\sigma_\delta(z)dz=1$, we define a smooth
function $\phi_\delta\in C^\infty_0(\mathbb{R}^{2n}, \mathbb{R}^+)$
by
\begin{equation}
\phi_\delta(z)=\frac{\sigma_\delta(z)}2+\frac{\sigma_\delta(N_0z)}2.
\end{equation}
It is clear that $\phi_\delta$ has compact support  in
$B(\delta/2)$, $\phi_\delta(N_0z)=\phi_\delta(z)$ and
$\int_{\mathbb{R}^{2n}}\phi_\delta(z)dz=1$.
 Define  a smooth function
$\rho:\mathbb{R}^{2n}\rightarrow \mathbb{R}$ with compact support in
 $A(2\delta)\subset Z(1)$ by
\begin{equation}
\rho(z)=\chi_{A(\delta)}*\phi_\delta(z)=\int_{\mathbb{R}^{2n}}\chi_{A(\delta)}(z-\tilde{z})\phi_\delta(\tilde{z})d\tilde{z}.
 \end{equation} From the definitions
we have $N_0A(\delta)=A(\delta)$,
$\chi_{A(\delta)}(N_0z)=\chi_{A(\delta)}(z)$,
$\phi_\delta(N_0z)=\phi_\delta(z)$, $\rho(N_0z)=\rho(z)$ and
$\rho|_{A(\delta/2)}\equiv 1$. Now  the Hamiltonian function
$K:Z(1)\rightarrow \mathbb{R}$  is defined by
\begin{equation}\label{K}
K(z)=\rho(z)\langle z, -J z_0\rangle.
\end{equation}
It is clear that $K(N_0z)=-K(z)$, $\nabla K(N_0z)=-N_0\nabla K(z)$
and $\nabla K(z)=-Jz_0$ for $z\in A(\delta/2)$, so
\begin{equation}\label{X_K}
X_K(N_0z)=J\nabla K(N_0z)=-JN_0\nabla K(z)=N_0J\nabla K(z)=N_0
X_K(z).
\end{equation}
The  flow $\psi^t:Z(1)\to Z(1)$ of the Hamitonian vector field $X_K$
is compact supported symplectic diffeomorphim for every $t>0$. We
define $\psi=\psi^1$ the time-$1$ map. It has compact support in
$Z(1)$, $\psi(N_0z)=N_0\psi(z)$,  $\psi^t (z)=z+tz_0$ for $z\in
A(\delta/2)$, so $\psi(0)=z_0$ and the Hamiltonian $H\cdot \psi$
vanishes in a neighborhood of $0$. The proof is complete.
$\hfill\Box$

From Lemma \ref{L2.6}, we only need to prove that the vector field
$X_{H\circ\psi}$ possesses a nonconstant  brake orbit with minimal
period $0<T\leq 1$.
 Hence we can assume that Hamiltonian $H$ vanishes in an open neighborhood of the
 origin.
We extend the function $H\in \mathcal {H}(Z(1))$ to a function
defined on the whole space $\mathbb{R}^{2n}$. This is possible since
$H$ is constant near the boundary of $Z(1)$.
 Denote by
\begin{equation}\label{q}
q(z)=q_K(z)=(x^2_1+y^2_1)+\frac{1}{K^2}\sum^n_{j=2}(x^2_j+y^2_j),
\end{equation}
where $z=(x,y)\in \mathbb{R}^{2n}$, and $K\in\mathbb{Z}^+$ is
sufficiently
 large.  It is clear that $q_K(z)=q_K(N_0z)$.
 Since $H\in\mathcal {H}(Z(1))$ there exists $K>0$ such that $H\in\mathcal {H}(E_K)$,
 where $E_K$ is defined by
$$E_K=\left\{z\in\mathbb{R}^{2n}|q_K(z)<1\right\}.$$
Since $H\in \mathcal{H}(Z(1))$ satisfies $m(H)>\pi$, there is an
$\varepsilon>0$ such that $m(H)>\pi+\varepsilon$. We can take a
smooth function $f:\mathbb{R}\rightarrow\mathbb{R}$ such that
$$f(s)=m(H)\;\;\;\; \text{for}\; s\leq 1,$$
$$f(s)\geq(\pi+\varepsilon)s\;\;\;\; \text{for all} \;s\in \mathbb{R},$$
$$f(s)=(\pi+\varepsilon)s\;\;\;\;  \text{for} \;s \;\text{large},$$
$$0<f'(s)\leq(\pi+\varepsilon)\;\;\;\; \text{for}\; s>1.$$
The extension of $H$ is now defined by
\begin{equation}\label{-H}
\bar{H}(z)=\left\{\begin{array}{ll} \;H(z), &z\in E_K,\\
 f(q_K(z)),&z\notin E_K.
 \end{array}\right.
\end{equation}
Clearly $\bar{H}\in C^{2}(\mathbb{R}^{2n})$ ,
$\bar{H}(N_0z)=\bar{H}(z)$ and $\bar{H}$ is quadratic at infinity,
exactly we have
\begin{equation}\label{-H2}
\bar{H}(z)=(\pi+\varepsilon)q_K(z),\;\;{\rm if }\;\; |z|\geq R
\end{equation}
for some large $R$.
 The following crucial lemma describes the
distinguished brake orbit we are looking for.
\begin{lem}\label{proposition 2}
Assume $x(t)$ is a brake orbit of
\begin{equation}\label{2.10}
\dot x(t)=J\bar{H}'(x(t))
\end{equation}
with period 1. If it satisfies
 \begin{equation}\label{2.11}
\Phi(x):=\int^1_0\left\{\frac{1}{2}\langle-J\dot
x(t),x(t)\rangle-\bar{H}(x(t))\right\}dt>0,
\end{equation}
 then $x(t)$ is
nonconstant and $x(t)\in E_K,\;\forall\,t\in [0,1]$. Hence $x(t)$ is
a nonconstant 1-periodic brake orbit of the original system $\dot
x(t)=JH'(x(t))$ on $Z(1)$.
\end{lem}
Since a brake orbit is a special periodic orbit, the proof of the
lemma is the same as the proof of Proposition 2 in [$P_{74}$,
\cite{H-Z}].

The remaining  of this subsection is to find a 1-periodic brake
orbit $x(t)$ of the equation \eqref{2.10} satisfying \eqref{2.11}.
We simply  replace $\bar{H}$ by $H$ in the sequel.

Denote by
\begin{equation}\label{L}
L^2=\{x\in L^2(S^1)|x=\sum_{j\in \mathbb{Z}}e^{2\pi jJt}x_j,x_j\in
\mathbb{R}^{2n},\;\sum_{j\in \mathbb{Z}}|x_j|^2<\infty\}.
\end{equation}
The space $ L^2$ is a Hilbert space with the usual $L^2$ inner
product $\langle\cdot,\cdot\rangle_0$  and associated norm
$\|\cdot\|_0$.
Denote by
\begin{equation}\label{Hs}
\widetilde{H}^s=\{x\in L^2(S^1)|x=\sum_{j\in \mathbb{Z}}e^{2\pi
jJt}x_j,x_j\in L_0,\;\sum_{j\in \mathbb{Z}}|j|^{2s}|x_j|^2<\infty\}.
\end{equation}
 The space $\widetilde{H}^s $ is a Hilbert space with
inner product and associated norm defined by
\begin{equation}\label{iner}
\langle x,y\rangle_s=\langle
x_0,y_0\rangle+2\pi\sum_{k\in\mathbb{Z}}|k|^{2s}\langle
x_k,y_k\rangle,
\end{equation}
\begin{equation}\label{norm}
\|x\|_s^2=\langle x,x\rangle_s,
\end{equation}
for $x,y\in \widetilde{H}^s$. Denote by $X=\widetilde{H}^{1/2}$,
$\|\cdot\|=\|\cdot\|_{1/2},\langle\cdot,\cdot\rangle=\langle\cdot,\cdot\rangle_{1/2}.$\\
\noindent There is an orthogonal splitting of $X$
\begin{equation}\label{ortho}
X=X^-\oplus X^0\oplus X^+
\end{equation}
with $X^+=\{x\in X|x=\displaystyle\sum_{j>0}e^{2\pi jJt}x_j\}$,
$X^-=\{x\in X|x=\displaystyle\sum_{j<0}e^{2\pi jJt}x_j\}$ and
$X^0=L_0$. The corresponding orthogonal projections are denoted by
$P^+,P^-,P^0$. Therefore, every $x\in X$ has a unique decomposition
$$x=x^-+x^0+x^+.$$
We define for $x,y\in X$
\begin{equation}\label{a(x,y)}
\begin{array}{ll}\vspace{0.5cm}
a(x,y)&=\frac{1}{2}\langle x^+,y^+\rangle_{1/2}-\frac{1}{2}\langle
x^-,y^-\rangle_{1/2}\\
&=\frac{1}{2}\langle(P^+-P^-)x,y\rangle_{1/2} ,
\end{array}
\end{equation}
which is a continuous bilinear form on $ X$. The functional $a:
X\rightarrow\mathbb{R}$, defined by
\begin{equation}\label{a(x)}
a(x)=a(x,x)=\frac{1}{2}\|x^+\|_{1/2}^2-\frac{1}{2}\|x^-\|_{1/2}^2,
\end{equation}
is differentiable with derivative

\begin{equation}\label{da}
da(x)(y)=\langle(P^+-P^-)x,y\rangle_{1/2},
\end{equation}
so  the gradient of $a$ is
\begin{equation}\label{a'}
\nabla a(x)=(P^+-P^-)x=x^+-x^-\in  X,\;\;\forall x\in  X.
\end{equation}
  We have $ X\subset  L^2$, the inclusion map
\begin{equation}\label{j}
j: X\rightarrow L^2
\end{equation}
is compact. Its adjoint operator
\begin{equation}\label{j*}
j^*:L^2\rightarrow X
\end{equation}
is defined by
\begin{equation}\label{jiner}
\langle j(x),y\rangle_0=\langle x,j^*(y)\rangle_{1/2},\;\;\forall
x\in  X,y\in L^2.
\end{equation}

\begin{lem}\label{Proposition 5} $j^*$ is
compact and there hold
$$j^*(L^2)\subset\widetilde{H}^1 \;\;\;\text{and}\;\;\; \|j^*(y)\|_1\leq\|y\|_0.$$
\end{lem}
\noindent{\it Proof.} By direct computation, we have for any $
y=\displaystyle\sum_{j\in\mathbb{Z}}e^{2\pi jJt}y_j\in L^2$,
$$j^*(y)=i^*(y_0)+\sum_{k\neq 0}\frac{1}{2\pi|k|}e^{k2\pi Jt}i^*(y_k),$$
where $i^*$ is the projection map: $\mathbb{R}^{2n}\rightarrow
L_0=\{0\}\oplus\mathbb{R}^n$. From the definition of $L^2$ and
$\widetilde{H}^s$, we can complete the proof.$\hfill\Box$

We next consider the functional
\begin{equation}\label{b}
b(x)=\int_0^1H(x(t))dt,
\end{equation}
since $H$ vanishes in a neighborhood of the origin, and from
\eqref{-H2},
  there is  $M>0$ such that
 $$|H''(z)|\leq M\;,\; |H'(z)|\leq M|z|\; {\rm and}\; |H(z)|\leq \frac{1}{2}M|z|^2, \;\;\forall z\in\mathbb{R}^{2n},$$
 so the functional $b$
 can be defined for $x\in L^2$ and hence also for $x\in  X\subset L^2$.

\begin{lem}\label{Lemma 4}
There holds
 $b\in C^1(X, \mathbb{R})$, $\nabla b: X\rightarrow X$  maps bounded sets into relatively compact sets.
 Moreover,
 $$\|\nabla b(x)-\nabla b(y)\|\leq M\|x-y\|$$
 and $|b(x)|\leq M\|x\|^2_{L^2},\forall x,y\in  X.$
\end{lem}
\noindent {\it Proof.} We have
 \begin{equation}\label{b'}
\nabla b(x)=j^*\nabla H(x).
 \end{equation}
 Moreover,
$$\begin{array}{ll} \|\nabla b(x)-\nabla b(y)\|_{1/2}&=\|j^*(\nabla
H(x)-\nabla
H(y))\|_{1/2}\\
\\
 &\leq\|\nabla H(x)-\nabla H(y)\|_{L^2}\\
 \\
 &\leq M\|x-y\|_{L^2}\\
 \\
&\leq M\|x-y\|_{1/2}.
\end {array}$$
 The proof is complete.$\hfill\Box$\\
 Now we consider the functional
 \begin{equation}\label{Phi}
\Phi(x)=a(x)-b(x),\; x\in X.
 \end{equation}
 We have $\Phi: X\rightarrow\mathbb{R}$ is differentiable
 and its gradient is given by
\begin{equation}\label{Phi'}
\nabla \Phi(x)=x^+-x^--\nabla b(x).
 \end{equation}

 \begin{lem}\label{Lemma 5}
Assume $x\in X$ is a
 critical point, i.e., $\nabla \Phi(x)=0$. Then $x\in C^2(S^1)$ and it
 is a brake orbit with 1-periodic.
 \end{lem}
 \noindent{\it Proof.} Let
 $$x(t)=\sum_{j\in \mathbb{Z}}e^{2\pi jJt}x_j,\;\;x_j\in \{0\}\oplus
 \mathbb{R}^n$$ be a critical point of $\Phi$, and
 $$\nabla H(x(t))=\sum_{j\in \mathbb{Z}}e^{j 2\pi Jt}a_j\in L^2,\;\;a_j\in
 \mathbb{R}^{2n}.$$
 So we have $$(P^+-P^-)x=j^*(\nabla H(x)).$$
 That is
\begin{equation}\label{x}
\left\{ \begin {array}{lll} i^*(a_0)=0,&k=0, \\
\frac{1}{2k\pi}i^*(a_k)=x_k,\;\;&k>0,
\\ -\frac{1}{2k\pi}i^*(a_k)=-x_k,\;\;&k<0,  \end{array}\right.
\end{equation}
 where $i^*$ is the projection map $\mathbb{R}^{2n}\to \{0\}\times \mathbb{R}^n$, so  $x\in \widetilde{H}^1$,
  $x(-t)=N_0x(t)$. Since $H(z)=H(N_0z)$, we have
$$N_0\nabla H(x(t))=\nabla H(N_0x(t))=\nabla H(x(-t)),$$
that is
$$N_0\sum_{j\in \mathbb{Z}}e^{2\pi jJt}a_j=\sum_{j\in \mathbb{Z}}e^{-2\pi jJt}a_j,$$
hence  $a_k\in \{0\}\oplus\mathbb{R}^n$, $a_k=i^*(a_k),\; \forall
k\in \mathbb{Z}.$ So from \eqref{x} we have $\dot x(t)=J\nabla
H(x(t))$, and  $x(t)$ is a brake orbit with 1-periodic.$\hfill\Box$
 \begin{lem}\label{Lemma 6}
 $\Phi$
 satisfies the (PS) condition.
 \end{lem}
\noindent{\it  Proof.} In fact we will prove that every sequence
 $\{x_j\}\subset X$ satisfying $\nabla \Phi(x_j)\rightarrow
 0$ contains a convergent subsequence. Assume $\nabla
\Phi(x_j)\rightarrow 0$ so that
\begin{equation}\label{}
x^+_j-x^-_j-\nabla b(x_j)\rightarrow 0.
\end{equation}
If $x_j$ is bounded in $ X$, then $x^0_j\in \mathbb{R}^{2n}$ is
bounded, and from Lemma \ref{Lemma 4}, we see that $\{x_j\}$ has
 a convergent subsequence. To prove that $x_j$ is
bounded we argue by contradiction and assume
$\|x_j\|\rightarrow\infty$. Define
\begin{equation}\label{}
y_k=\frac{x_k}{\|x_k\|},
\end{equation}
so
 $\|y_k\|=1$. By assumption, from \eqref{b'},
$$(P^+-P^-)y_k-j^*(\frac{1}{\|x_k\|}\nabla H(x_k))\rightarrow 0.$$
Since $|\nabla H(z)|\leq M|z|$, the sequence
$$\frac{\nabla H(x_k)}{\|x_k\|}\in L^2$$
is bounded in $L^2$. Since $j^*:L^2\rightarrow X$ is compact,
$(P^+-P^-)y_k$ is relatively compact, and since $y^0_k$ is bounded
in $\mathbb{R}^{2n}$, the sequence $y_k$ is relatively compact in $
X$. After taking a subsequence we can assume $y_k\rightarrow y$ in $
X$ and hence $y_k\rightarrow y$ in $L^2$. From \eqref{-H2}, we have
$$\left\|\frac{\nabla  H(x_k)}{\|x_k\|}-\nabla Q(y)\right\|_{L^2}\leq\frac{1}{\|x_k\|}\|\nabla H(x_k)-\nabla Q(x_k)\|_{L^2}+\|\nabla Q(y_k-y)\|_{L^2},$$
where $Q(z)=(\pi+\varepsilon)q(z)$. Since $|\nabla H(z)-\nabla
Q(z)|\leq M$ for all $z\in \mathbb{R}^{2n}$ and since $\nabla Q$
defines a continuous linear operator of $L^2$, we conclude
$$\frac{\nabla H(x_k)}{\|x_k\|}\rightarrow \nabla Q(y) \;\;\text{in}\;\; L^2.$$
Consequently,
$$\frac{\nabla b(x_k)}{\|x_k\|}=j^*(\frac{\nabla H(x_k)}{\|x_k\|})\rightarrow j^*(\nabla Q(y))\;\;\text{in}\;\;  X.$$
This implies that $y\in  X$ solves the linear equation in $ X$
$$y^+-y^--j^*\nabla Q(y)=0,$$
$$\|y\|=1.$$
As in Lemma \ref{Lemma 5} one verifies that $y$  solves the linear
Hamiltonian equation
$$\dot y(t)=J\nabla Q(y(t)).$$
Recall now that $Q=(\pi+\varepsilon)q$, and
$q(z)=(x^2_1+y^2_1)+\frac{1}{K^2}\sum^n_{j=2}(x^2_j+y^2_j)$. We see
that the symplectic 2-planes $\{x_j,y_j\}$ are filled with periodic
solutions of $J\nabla Q$ having periods $T\neq 1$. Since the linear
equation does not admit any nontrivial periodic solutions of period
1 we conclude $y(t)\equiv 0$. This contradicts $\|y\|=1$ and we
conclude that the sequence $\{x_k\}$ must be bounded.$\hfill\Box$

$\nabla \Phi$ is  globally Lipschitz continuous, so the  gradient
equation
$$\dot x=-\nabla \Phi(x),\;\;x\in X$$
 defines a unique
global flow
$$\mathbb{R}\times X\rightarrow X \;\; : \;\; (t,x)\mapsto\varphi^t(x)\equiv x\cdot t,$$
which maps bounded sets into bounded sets.
\begin{lem}\label{Lemma 7}
The flow of $\dot x=-\nabla \Phi(x)$ has the following form
\begin{equation}\label{x-K}
x\cdot t=e^tx^-+x^0+e^{-t}x^++K(t,x),
\end{equation}
where $K:\mathbb{R}\times  X\rightarrow  X$ is continuous and maps
bounded sets into precompact sets.
\end{lem}
\noindent{\it Proof.} Define a map $K$ by
\begin{equation}\label{K}
K(t,x)=\int^t_0(e^{t-s}P^-+P^0+e^{-t+s}P^+)\nabla b(x\cdot s)ds.
\end{equation}
We have to verify that $K$ has the desired properties. Denote the
right hand side of \eqref{x-K} by $y(t)$, we have that
$$\dot y(t)=(P^--P^+)y(t)+\nabla b(x\cdot t).$$
Since $y(0)=x$, the function $\xi(t)=y(t)-x\cdot t$ solves the
linear equation
$$\dot\xi(t)=(P^--P^+)\xi(t) \;\;\text{and}\; \xi(0)=0.$$
By the uniqueness of the initial value problem $\xi(t)=0$ so that
$y(t)=x\cdot t$ as required. In view of \eqref{b'} we can write
$$K(t,x)=j^*\{\int^t_0(e^{t-s}P^-+P^0+e^{-t+s}P^+)\nabla H(j(x\cdot s))ds\}.$$
 By Lemma\ref{Proposition 5},
$j^*:L^2\rightarrow X$ maps bounded sets into precompact sets and,
therefore, $K$ has the desired properties.$\hfill\Box$
\begin{prop}\label{Proposition 6}
There exists $x^*\in  X$ satisfying $\nabla \Phi(x^*)=0$ and
$\Phi(x^*)>0$
\end{prop}
In order to prove this proposition we first single out two
 subsets $\Omega$ and $\Gamma$ of $ X$. The bounded set
$\Omega=\Omega_\tau\subset X$ is defined by
\begin{equation}\label{Sigma}
\Omega_\tau=\{x|x=x^-+x^0+se^+,\|x^-+x^0\|\leq
\tau\;\;\text{and}\;\; 0\leq s\leq\tau\},
\end{equation}
where $\tau>0$ and $e^+\in X^+$ is defined by
$$e^+(t)=e^{2\pi Jt}e_1\;\;\text{and}\;\;e_1=(0,...,0,1,...,0)\in \{0\}\oplus\mathbb{R}^{n}.$$
Clearly $\|e^+\|^2=2\pi$ and $\|e^+\|_{L^2}=1$.   We denote
$\partial\Omega$ the boundary of $\Omega$ in $X^-\oplus
X^0\oplus\mathbb{R}e^+$.
\begin{lem}\label{Lemma 8}
There exists $\tau^*>0$ such that for $\tau>\tau^*$
$$\Phi|_{\partial\Omega_{\tau}}\leq 0.$$
\end{lem}
\noindent{\it Proof.} From $a|_{X^-\oplus X^0}\leq 0$ and $b\geq 0$
we have
$$\Phi|_{X^-\oplus X^0}\leq 0.$$
We shall deal with the functional on those parts of the boundary
$\partial\Omega_{\tau}$ which are defined by $\|x^-+x^0\|=\tau$ or
$s=\tau$. By the construction of $H$ there exists a constant
$\gamma>0$ such that
$$H(z)\geq (\pi+\varepsilon)q(z)-\gamma\;\;\text{ for all}\;\;z\in\mathbb{R}^{2n}.$$
Therefore,
$$\Phi(x)\leq a(x)-(\pi+\varepsilon)\int^1_0q(x)+\gamma,\;\;\text{for all}\;\; x\in X.$$
Recalling the definition of the quadratic form $q$, one verifies for
$x=x^-+x^0+se^+\in X^-\oplus X^0\oplus X^+$ that
$$\int^1_0q(x^-+x^0+se^+)dt=\int^1_0q(x^-)dt+\int^1_0q(x^0)dt+\int^1_0q(se^+)dt.$$
Recalling that $\|e^+\|^2=2\pi$ , for $x=x^-+x^0+se^+$, there holds
$$\begin{array}{ll}
\Phi(x)&\leq\frac{1}{2}s^2\|e^+\|^2-\frac{1}{2}\|x^-\|^2-(\pi+\varepsilon)q(x^0)-(\pi+\varepsilon)\int^1_0q(se^+)+\gamma\\
\\
 &=-\frac{1}{2}\|x^-\|^2-\varepsilon
 s^2\|e^+\|^2_{L^2}-(\pi+\varepsilon)q(x^0)+\gamma.
\end {array}$$
Consequently there exists a constant $c>0$ such that
$$\Phi(x^-+x^0+se^+)\leq\gamma-c\|x^-+x^0\|^2-c\|se^+\|.$$
The right hand side is not positive if $\|x^-+x^0\|=\tau$ or
$s=\tau$ for $\tau$ sufficiently large. The proof of the lemma
is complete.$\hfill\Box$

 The subset $\Gamma=\Gamma_\alpha\subset
X^+$ is defined by
\begin{equation}\label{Gamma}
\Gamma_{\alpha}=\{x\in X^+\mid\|x\|=\alpha\}.
\end{equation}
\begin{lem}\label{Lemma 9}
There exist  $\alpha>0$ and $\beta>0$ such that
$$\Phi|_{\Gamma_{\alpha}}\geq\beta>0.$$
\end{lem}
\noindent{\it Proof.}  The space $ X$ is continuously embedded in
$L^p(S^1)$ for every $p\geq 1$. Hence there is a constant $M=M_p$
such that
$$\|u\|_{L^p}\leq M\|u\|_{1/2},\;\;\;u\in  X.$$
Observing that $|H(z)|\leq c|z|^3$ for all $z\in\mathbb{R}^{2n}$, we
can take a constant $K>0$ such that
$$\int^1_0|H(x(t))|dt\leq c\|x\|^3_{L^3}\leq K\|x\|^3_{1/2},$$
for all $x\in  X$. Now, if $x\in X^+$, then $\Phi(x)\geq
\frac{1}{2}\|x\|^2-K\|x\|^3$ and the lemma is now obvious for some
small $\alpha>0$ and $\beta>0$.$\hfill\Box$

 Since $\Phi(\varphi^t(x))$ decreases in $t$ we
conclude immediately from Lemma \ref{Lemma 8} and Lemma \ref{Lemma
9} that $\varphi^t(\partial\Omega)\cap\Gamma=\emptyset$ for all
$t\geq 0$.  But the following result tell us that
$\varphi^t(\Omega)\cap\Gamma\ne\emptyset$.

\begin{lem}\label{Lemma 10}
$$\varphi^t(\Omega)\cap\Gamma\neq\emptyset,\;\; \forall\;\;t\geq 0.$$
\end{lem}
 \noindent {\it Proof.} We shall use the
 Leray-Schauder degree. Abbreviating the flow by $\varphi^t(x)\equiv x\cdot
 t$, we need to verify that $(\Omega\cdot t)\cap\Gamma\neq
 \emptyset$ for all $t\geq 0$. We can rewrite this by requiring
 \begin{equation}\label{x-in-sigma}
 (P^-+P^0)(x\cdot t)=0,\;\;\;
 \|x\cdot t\|=\alpha,\;\;\;
 x\in\Omega.\;
 \end{equation}
 Recall that, by Lemma \ref{Lemma 7}, the flow has the
 representation $x\cdot t=e^tx^-+x^0+e^{-t}x^++K(t,x)$, so that
 \eqref{x-in-sigma} becomes
\begin{equation}\label{x-in-sigma-2}
 \begin{array}{rrr} e^tx^-+x^0+(P^-+P^0)K(t,x)=0,\;\;\;
 \alpha-\|x\cdot t\|=0,\;\;\;
 x\in\Omega. \end {array}
 \end{equation}
 Multiplying the $X^-$ part by $e^{-t}$ one gets the following equivalent
 equations
\begin{equation}\label{x-in-sigma-3}
  x^-+x^0+(e^{-t}P^-+P^0)K(t,x)=0,\;\;\;
 \alpha-\|x\cdot t\|=0,\;\;\;
 x\in\Omega.\;\;\;
 \end{equation}
 Since $x\in\Omega$ is represented by $x=x^-+x^0+se^+$, with $0\leq
 s\leq\tau$, we can rewrite \eqref{x-in-sigma-3} as follows:
 \begin{equation}\label{x-B-sigma}
x+B(t,x)=0\;\;\text{and}\;\; x\in\Omega,
 \end{equation}
where the operator $B$ is defined by
\begin{equation}
B(t,x)=(e^{-t}P^-+P^0)K(t,x)+P^+\{(\|x\cdot t\|-\alpha)e^+-x\}.
\end{equation}
Abbreviating $F=X^-\oplus X^0\oplus\mathbb{R}e^+$, the map
$B:\mathbb{R}\times F\rightarrow F$ is continuous and maps bounded
sets into relatively compact sets. This was proved in Lemma
\ref{Lemma 7}. We therefore can apply the Leray-Schauder degree
theory. The equation \eqref{x-B-sigma}
 has a solution $x\in\Omega$ for given $t\geq 0$ if  deg$(\Omega,id+B(t,\cdot),0)\neq 0$. In view of
 $\varphi^t(\partial\Omega)\cap\Gamma=\emptyset$ for $t\ge 0$, we have
 \begin{equation}\label{}
0\notin(id+B(t,\cdot))(\partial\Omega),\;\forall\,t\geq 0.
 \end{equation}
 Hence by the homotopic invariance of the degree, there holds
\begin{equation}
\text{deg}(\Omega,id+B(t,\cdot),0)=\text{deg}(\Omega,id+B(0,\cdot),0).
\end{equation}
Since $K(0,x)=0$ we find $B(0,x)=P^+\{(\|x\|-\alpha)e^+-x\}$.
Defining the homotopy
\begin{equation}
L_\mu(x)=P^+\{(\mu\|x\|-\alpha)e^+-\mu x\}\;\; \text{for}\;\; 0\leq
\mu\leq 1,
\end{equation}
we claim $x+L_\mu(x)\neq 0$ for $x\in\partial\Omega$. Indeed, if
$x\in\Omega$ satisfies $x+L_\mu(x)=0$ then $x=se^+$ and, therefore,
$s((1-\mu)+\mu\|e^+\|)=\alpha$. Consequently $0<s\leq\alpha$, so  as
claimed $x\notin\partial \Omega$ if $\tau>\alpha$. Therefore, by
homotopic invariance again, there holds
$$
\text{deg}(\Omega,id+B(t,\cdot),0)=\text{deg}(\Omega,id+L_0,0)
=\text{deg}(\Omega,id-\alpha
B(t,\cdot)e^+,0)\\
=\text{deg}(\Omega,id, \alpha e^+)=1
$$
\noindent provided that $\alpha e^+\in\Omega$, which holds true for
$\tau>\alpha$. This finishes the proof of Lemma \ref{Lemma
10}.$\hfill\Box$

Now we can finish the proof of Proposition \ref{Proposition 6}. We
shall apply the  minimax argument. We take the family $\mathcal {F}$
consisting of the subsets $\varphi^t(\Omega)$, for every $t\geq 0$
and define
\begin{equation}\label{minimax}
c(\Phi,\mathcal {F})=\inf_{t\geq
0}\;\sup_{x\in\varphi^t(\Omega)}\Phi(x).
\end{equation}
We claim that $c(\Phi,\mathcal {F})$ is finite. Indeed, since
$\varphi^t(\Omega)\cap\Gamma\neq\emptyset$ and
$\Phi|_{\Gamma}\geq\beta$ we conclude that
\begin{equation}\label{ineq}
\beta\leq\inf_{x\in\Gamma}\Phi(x)\leq\sup_{x\in\varphi^t(\Omega)}\Phi(x)<\infty.
\end{equation}
In  the last estimate of \eqref{ineq} we have used that $\Phi$ maps,
in view of Lemma \ref{Lemma 4}, bounded sets into bounded sets.
Therefore,
\begin{equation}
-\infty<\beta\leq c(\Phi,\mathcal {F})<\infty.
\end{equation}
 We know already that the functional $\Phi$ satisfies the (PS)
 condition (Lemma\ref{Lemma 6}). Moreover,  the family $\mathcal{F}$ is  invariant under the negative gradient
 flow $\varphi^t$ for $t>0$.
  Consequently the Minimax Lemma implies that $c(\Phi,\mathcal {F})$ is a critical value. We
 deduce that there is a point $x^*\in X$ satisfying $\nabla\Phi(x^*)=0$
 and
 $$\Phi(x^*)=c(\Phi,\mathcal {F})\geq\beta>0,$$
 and the proof of Proposition \ref{Proposition 6} is complete. \\
\subsection{Application to the Existence of Brake Orbit}
\label{2.2}
 In this subsection, we  use the
symmetrical symplectic capacity theory developed in the previous
subsection to solve the existence
of brake orbits on energy surfaces.

 Let $(\mathcal{M},\mathcal{L},\omega, \varphi)$ be a symmetrical
symplectic manifold, and $H\in C^2(\mathcal{M},\mathbb{R})$
satisfying $H(\varphi(x))=H(x),\;\forall\; x\in \mathcal{M}$.
Suppose that the energy surface
\begin{equation}\label{S}
\Sigma=\{ x\in\mathcal{M}| H(x)=1 \}
\end{equation}
is compact and regular, i.e.,
 \begin{equation}\label{dH}
dH(x)\neq 0 \;\;{\rm for}\;\; x\in \Sigma,
 \end{equation}
 and $\Sigma\cap \mathcal {L}\neq \emptyset$ with transversal intersections. Thus $\Sigma\subset \mathcal{M}$ is a smooth and compact submanifold
 of codimension 1 whose tangent space at $x\in \Sigma$ is given by
 \begin{equation}\label{TxS}
T_x\Sigma=\{\xi\in T_x\mathcal{M}|dH(x)\xi=0 \}.
 \end{equation}
We define  an open and bounded neighborhood $U$ of $\Sigma$ by
\begin{equation}\label{U}
U=\bigcup_{\lambda\in I}\Sigma_\lambda,
\end{equation}
where $I=(1-\varepsilon,1+\varepsilon)$ for some small
$\varepsilon>0$, and  $\Sigma_\lambda=\{x\in \mathcal
{M}|H(x)=\lambda\}$ is diffeomorphic to  $\Sigma$ with
$\Sigma_{\lambda}\cap \mathcal {L}\neq \emptyset$ for all
$\lambda\in I$. Indeed, the gradient
 $\nabla H\neq 0$ in a neighborhood of
$\Sigma$, in view of \eqref{U}. The modified gradient flow
$\psi_0^t$ defined by the following equation
$$\dot{x}=\frac{\nabla H(x)}{|\nabla H(x)|^2}$$
 is transversal to $\Sigma$, and there holds
$$H(\psi_0^t(x))=1+t, \;\forall\,x\in \Sigma.$$
This means that $\psi_0^t:\Sigma\to \Sigma_{1+t}$ is a
diffeomorphism.  Since $H(\varphi(z))=H(z)$, we have $\varphi(U)=U$.
Similar to Theorem 1 in $P_{106}$ of \cite{H-Z}, we have the
following result which is equivalent to Theorem \ref{H-Z}.

\begin{thm}\label{c-H-Z}
There is a dense subset $O\subset I$, such that for $\lambda\in O$
the energy surface $\Sigma_\lambda$ possesses a brake orbit of
$X_H$, provided $c_{\varphi}(U,\omega)<\infty$.
\end{thm}
 \noindent{\it Proof.}
Suppose $I=(1-\rho,1+\rho)$ for some small $\rho>0$. For
$0<\varepsilon<\rho$, we define a smooth function
$f:\mathbb{R}\rightarrow\mathbb{R}$  by
$$\left\{\begin{array}{lll}
f(s)&=c_{\varphi}(U)+1,\;\;&{\rm for}\;\; s\leq 1-\varepsilon
\;\;{\rm and
}\;\;s\geq  1+\varepsilon,\\
\\
f(s)&=0,\;\;&{\rm for }\;\;1-\frac{\varepsilon}{2}\leq s\leq
1+\frac{\varepsilon}{2},\\
\\
f'(s)&<0,\;\;&{\rm for }\;\;1-\varepsilon<s<1-\frac{\varepsilon}{2},\\
\\
f'(s)&>0,\;\;&{\rm for }\;\;1+\frac{\varepsilon}{2}<s<1+\varepsilon.
\end{array}\right.$$

\noindent Define $F:U\rightarrow \mathbb{R}$ by
$$F(x)=f(H(x)),\;\;x\in U.$$
It is easy to see that $F\in \mathcal
{H}(U,U\cap\mathcal{L},\omega,\varphi)$ and $m(F)>c_{\varphi}(U)$.
Consequently, in view of the definition of the capacity
$c_{\varphi}(U)$, there exists a nonconstant brake orbit $(T,x(t))$
with $0<T\leq 1$ of the Hamiltonian system:
$$\dot{x}=X_F(x(t)),\;\;x(t)\in U,$$
where
$$X_F(x)=f'(H(x))\cdot X_H(x),\;\;x\in U.$$
Moreover,
$$H(x(t))=\lambda$$
is constant in $t$. Since $x(t)$ is not a constant solution we
conclude
$$f'(H(x(t)))=f'(\lambda)=\tau\neq 0.$$
Thus, in view of the definition of the function $f$, the value
$\lambda$ belongs to the set
$1-\varepsilon<\lambda<1-\frac{\varepsilon}{2}$ or
$1+\frac{\varepsilon}{2}<\lambda<1+\varepsilon$. In particular
$|\lambda-1|<\varepsilon$. By rescaling, we define
$$y(t)=x(\frac{t}{\tau}),$$
which has period $|\tau| T$ and satisfies
$$\dot y(t)=X_H(y(t)),$$
hence $y(t)$ is a brake orbit of the original Hamiltonian vector field
$X_H$ on the energy surface $H(y(t))=\lambda$. Duo to the
arbitrariness of  $\varepsilon$, $1$ is the limit point of $\lambda$
such that $\Sigma_{\lambda}$ possesses a brake orbit and so is true for all point in $I$. $\hfill\Box$ \\

\begin{rem}
Actually  the Lebesgue measure of $O$ in Theorem \ref{c-H-Z} is
equal to Lebesgue measure of $I$, i.e., $m(O)=m(I)$. The proof is
similar to that of Theorem 4 in $P_{118}$ of \cite{H-Z}.
\end{rem}

\begin{defi}\label{3.4}
We call a hypersurface $\Sigma\subset \mathcal{M}$ is
$\varphi$-invariant, if it satisfies $\varphi (\Sigma)=\Sigma$ and
$\Sigma\cap \mathcal{L}\neq\emptyset$. We denote by
$\mathcal{S}_\varphi$ the set of all $\varphi$-invariant
hypersurface in $\mathcal{M}$.
\end{defi}
\begin{defi}\label{r-c-t}
A compact  hypersurface $\Sigma\in \mathcal{S}_\varphi$ is called
$\varphi$-contact type if there exists a vector field $X$, defined
on a neighborhood $U$ of  $\Sigma$, and a constant $\lambda\neq 0$
such that
\begin{equation}\label{r-c-t-2}
\left\{\begin{array}{lll}
L_X\omega=\lambda\;\omega,\;\;on \;U,\\
\\
X(x)\notin T_x \Sigma,\;\;\forall\; x\in \Sigma,\\
\\
\varphi_* (X(x))=X(\varphi(x)),\;\; \forall\,x\in U.
\end{array}\right.
\end{equation}

\end{defi}

\noindent {\it Proof of Theorem \ref{C.Viterbo}.} We follow the
ideas of the proofs of Theorem 5 and Theorem 6 in ${\rm P}_{123}$ of
\cite{H-Z}. Let $X$ be the  vector field defined in Definition
\ref{r-c-t}. Since $\Sigma$ is compact and $X$ is transversal to
$\Sigma$, the map
\begin{equation}\label{X-flow}
\Psi:\Sigma\times(-\varepsilon,\varepsilon)\rightarrow U\subset
\mathcal{M}
\end{equation}
defined by $\Psi(x,t)=\psi^t(x)$ for $x\in \Sigma$ and
$\|t\|<\varepsilon$ is a diffeomorphism onto an open neighborhood
$U$ of $\Sigma$ provided $\varepsilon>0$ is sufficiently small,
where $\psi^t$ is the flow of $X$. From $L_X\omega=\omega$ we
conclude that if $x(s)$ is a closed brake characteristic on
$\Sigma$, then $y(s)=\psi^t(x(s))$ will be a closed brake
characteristic on $\Sigma_t=\psi^t(\Sigma)$, then from Theorem
\ref{c-H-Z}, we complete the proof.$\hfill\Box$

 Consider the Example
\ref{eg1.1}, from the results above,  we have the following result.

\begin{cor}\label{3.7}
Let $H\in C^2(\mathbb{R}^{2n},\mathbb{R})$ with
$H(N_0x)=H(x),\forall x\in\mathbb{R}^{2n}$. Suppose
$\Sigma=H^{-1}(1)$ is its  compact  regular $N_0$-invariant energy
surface, $\Sigma\cap \mathcal {L}\neq \emptyset$ with transversal
intersections. Then for an open interval
$I=(1-\varepsilon,1+\varepsilon)$, $\varepsilon>0$ small, there is a
dense subset $O\subset I$ such that for all $\lambda\in O$ the
energy surface $\Sigma_\lambda=H^{-1}(\lambda)$ possesses  a brake
orbit of $X_H$. Moreover if $\Sigma$ is $N_0$-contact type, then it
carries a  brake orbit.$\hfill\Box$
\end{cor}
\begin{rem}
 It is easy to say that if $\Sigma$ is $N_0$-invariant and star-shaped with center at origin, then $\Sigma$ is
 $N_0$-contact type, and Corollary \ref{3.7} generalize the result of
 \cite{R}.
\end{rem}
\noindent For further applications, we need the following lemma:
\begin{lem}\label{3.8}
A compact hypersurface $\Sigma\in \mathcal{S}_\varphi$ is of
$\varphi$-contact type if and only if there exists a 1-form $\alpha$
on a neighborhood $U$ of $\Sigma$  and  constant $\lambda\neq 0$
such that
\begin{equation}
\left\{ \begin{array}{lll}
d\alpha=\lambda\;\omega,\\
\alpha(\xi)\neq 0, \;\;{\rm for} \;\;0\neq\xi\in\mathcal{L}_\Sigma,\\
\varphi^*(\alpha)(x)=-\alpha(\varphi(x)), \;\;\forall z\in U,
\end{array} \right.
\end{equation}
 where $\mathcal{L}_\Sigma=\{(x,\xi)\in T
\Sigma| \omega_x(\xi,\eta)=0\,\;\;\forall \eta\in T_x\Sigma\}$.
\end{lem}
\noindent{\it Proof.} Let $\alpha=i_X\omega$, since
$\varphi^*\omega=-\omega$, we have
$$
\varphi^*(\alpha)=i_{\varphi_*X}(\varphi^*\omega)=-i_{\varphi_*X}\omega,
$$
that is to say
$$
\varphi^*(\alpha)=-\alpha\Leftrightarrow\varphi_*(X)=X.
$$
The  remains of the proof is similar to \cite{H-Z}.$\hfill\Box$

 \begin{lem}\label{3.9}
Let $H\in C^\infty(\mathbb{R}^{2n},\mathbb{R})$,
$H(N_0(x,y))=H(x,y)$ satisfying
\begin{equation}\label{311}
\langle\frac{\partial}{\partial x}H(x,y),x\rangle>0,\;\forall
(x,y)\in\mathbb{R}^{2n}, {\rm \;with\;} x\neq 0.
\end{equation}
  Then  every compact and regular energy surface
$\Sigma=H^{-1}(c)$ with
\begin{equation}\label{3.12}
c<  \sup_{y\in\mathbb{R}^n}H(0,y)
\end{equation}
 belongs to $\mathcal{S}_{N_0}$
(see Definition \ref{3.4}) and is of $N_0$-contact type.
 \end{lem}
 \noindent{\it Proof.} Since $H$ is $N_0$-invariant, we have
 $N_0(\Sigma)=\Sigma$. Let $(x_0,y_0)\in \Sigma$, that is $H(x_0,y_0)=c$, if $x_0=0$,
  we have  $\Sigma\cap L_0\neq\emptyset$.
  Otherwise  from
 \eqref{311} we have $H(0,y_0)<c$, and from \eqref{3.12},
 there exists a $(0,y_1)\in L_0$ such that $H(0,y_1)>c$. So from the
 smoothness of $H$,  there exists a point $(0,y)\in L_0$ such that $H(0,y)=c$,
 so we also have $\Sigma\cap L_0\neq\emptyset$.
  In order to show that $\Sigma$ is $N_0$-contact type, we define the
  1-form on $\mathbb{R}^{2n}$ as in Lemma \ref{3.8} by
\begin{equation}
\alpha_\varepsilon=-xdy+\varepsilon dF\in
T^*_{(x,y)}\mathbb{R}^{2n},\;\varepsilon\in\mathbb{R},
\end{equation}
where $F\in C^\infty(\mathbb{R}^{2n},\mathbb{R})$ defined by
\begin{equation}
F(x,y)=\langle x,\frac{\partial}{\partial y}H(0,y)\rangle.
\end{equation}
From the definition of $\alpha_\varepsilon$ we have
\begin{equation}\label{3.16}
d\alpha_\varepsilon=\omega_0.
\end{equation}
Clearly $F(N_0(x,y))=-F(x,y)$ and¡¡$\;N_0^*dF=-dF$, so
\begin{equation}\label{3.17}
N_0^*\alpha_\varepsilon=-\alpha_\varepsilon.
\end{equation}
Finally since $\Sigma$ is a regular energy surface of $H$, and
$H(N_0(x,y))=H(x,y)$, we have if $(0,y)\in \Sigma$,
\begin{equation}
X_H(0,y)=J\nabla H(0,y)=(-\frac{\partial}{\partial
y}H(0,y),\frac{\partial}{\partial
x}H(0,y))^T=(-\frac{\partial}{\partial y}H(0,y),0)^T,
\end{equation}
that is to say
\begin{equation}\label{3.19}
\frac{\partial}{\partial y}H(0,y)\neq 0, {\rm \;if\;} (0,y)\in
\Sigma. \end{equation} By definition,  we have
\begin{equation}
\alpha_\varepsilon(X_H)=-\langle\frac{\partial}{\partial
x}H(x,y),x\rangle-\varepsilon\langle\frac{\partial}{\partial
y}H(x,y),\frac{\partial}{\partial y}H(0,y)\rangle+\varepsilon
\sum_{j,k=1}^n\frac{\partial^2}{\partial y_k\partial y_j
}H(0,y)\frac{\partial}{\partial x_k}H(x,y)x_j.
\end{equation}
So from \eqref{311}, \eqref{3.19} and the compactness of $\Sigma$,
 there is a $\delta>0$ small enough, such that
\begin{equation}\label{3.21}
\alpha_\varepsilon(X_H)(x,y)<0, {\rm \;if\; }(x,y)\in \Sigma, {\rm
\; with\; }\|x\|<\delta,
\end{equation}
and also from \eqref{311} and the compactness of $\Sigma$, there
exists a $\varepsilon$ depending on the $\delta$, such that
\begin{equation}\label{3.22}
\alpha_\varepsilon(X_H)(x,y)<0, {\rm \;if\; }(x,y)\in \Sigma, {\rm
\; with\; }\|x\|\geq\delta.
\end{equation}
Then from \eqref{3.16}, \eqref{3.17}, \eqref{3.21} and \eqref{3.22},
one see that $\Sigma$ is of $N_0$-contact type.$\hfill\Box$

\noindent So from Corollary \ref{3.7} and Lemma \ref{3.8},
\ref{3.9}, we have the following result.
\begin{thm}\label{3.11}
Let $H\in C^\infty(\mathbb{R}^{2n},\mathbb{R})$, $H(-x,y)=H(x,y)$
satisfying
\begin{equation}\label{3.20}
\langle\frac{\partial}{\partial x}H(x,y),x\rangle>0,\;\forall
(x,y)\in\mathbb{R}^{2n}, {\rm \;with\;} x\neq 0.
\end{equation}
  Then  every compact  regular energy surface
$\Sigma=H^{-1}(c)$ with
$$ c<  \sup_{y\in\mathbb{R}^n}H(0,y) $$
possesses  a  brake orbit of $X_H$.
 \end{thm}

For the $\varphi$-symmetric symplectic manifold
$(T^*(T^n),\mathcal{L},\omega,\varphi)$ discussed in Example
\ref{eg1.3} with its coordinates  $(x,y)$,  $x\in T^n$ and $y\in
T^*_x(T^n)\simeq\mathbb{R}^{n}$. We have the following result.

\begin{thm}  Suppose $H\in
C^2(T^*(T^n),\mathbb{R})$ satisfying $H(-x,y)=H(x,y)$ and $H(x,y)\to
+\infty$ with $|y|\to \infty$. Then for a
 regular energy
hypersurface $\Sigma_s=H^{-1}(s)$ with $s>\displaystyle\min_{y\in
\mathbb{R}^{n}}H(0,y)$, there is a sequence $s_k\to s\;(k\to
+\infty)$ such that the energy surface $\Sigma_{s_k}=H^{-1}(s_k)$
possesses  at least one brake orbit.
\end{thm}

\noindent{\it Proof.} By the condition $H(x,y)\to +\infty$ ($|y|\to
\infty$), we see that the energy hypersurface $H^{-1}(s)$ is compact
and
 there is a constant $a>0$ such that $\displaystyle\bigcup_{t\in[s-\delta,s+\delta]} H^{-1}(t)\subseteq
 T^n\times(-a,a)^n$ for some  $\delta>0$. Combining with $s>\displaystyle\min_{y\in
\mathbb{R}^{n}}H(0,y)$ we have $\mathcal{L}\cap \Sigma_s\neq
\emptyset$. So the result of this theorem comes from Theorem
\ref{H-Z} and the following result. $\hfill\Box$

\begin{lem}\label{Tn}
Let $(T^*(T^n),\mathcal{L,\omega,\varphi})$ be a symmetrical
symplectic manifold defined in Example \ref{eg1.3}, then we have
$$c_\varphi(T^n\times(-a,a)^n)\leq 5a\pi, \forall a>0.$$
\end{lem}
\noindent {\it Proof.} In fact, we can get a symplectic
diffeomorphism $\phi$ from \cite{H-Z} and \cite {Jiang}
$$\phi: S^1\times (-a,a)\rightarrow A=\left\{(x,y)\in
\mathbb{R}^2|a<x^2+y^2<5a\right\},$$ by
$$\phi(\theta,r)=((3a+r)^{1/2}cos\theta,(3a+r)^{1/2}sin\theta),$$
where  $0\leq \theta\leq2\pi$ and $-a<r<a$. Then it is easy to
verify that $\phi^*(dy\wedge dx)=dr\wedge d\theta$. Extending to
high dimensional case in the obvious way, we get a symplectic
diffeomorphism
$$\Phi:T^n\times(-a,a)^n\rightarrow A\times A\times\cdots\times A\subset \mathbb{R}^{2n}$$
satisfying $\Phi\varphi=N_1\Phi$, where $N_1=\left(\begin{matrix}I_n
&0\\0
 &-I_n\end{matrix}\right)$.
 So from Theorem \ref{main-result}, we have
$$c_\varphi(T^n\times(-a,a)^n,\omega)=c_{N_1}(A\times A\times\cdots\times A,\omega_0)\leq c_{N_1}(B((5a)^{1/2}),\omega_0),$$
and  $c_{N_1}(B((5a)^{1/2}),\omega_0)=5a
c_{N_1}(B(1),\omega_0)=5a\pi$.$\hfill\Box$

\section{($N_0,S$)-Symmetrical Symplectic Capacity and Applications}
In this section, we consider the following problem on $U$,
\begin{equation}\label{s-brake-orbit}
 \left\{ \begin {array}{lll}
\dot{x}(t)=J\nabla H(x(t)),\\
 x(-t)=N_0x(t),\\
x(t+\frac{T}{m})=Sx(t),
  \end{array}\right.
 \end{equation}
 where  $J$ and $N_0$ are
 defined in Example \ref{eg1.1}, and $S$ is an orthogonal symplectic
  matrix satisfying $S^m =I_{2n\times 2n}$ with a fixed
  constant $m\in \mathbb{N}\backslash\{1\}$.
  $U$ is an open subset of $(\mathbb{R}^{2n}, \omega_0)$, satisfying $\{0\}\in U$,
  $N_0U=SU=U$ ($(N_0,S)$-invariant subset).
  $H\in C^2(U,\mathbb{R})$,
  satisfying $H(N_0z)=H(Sz)=H(z),\;\forall z\in U$ ($(N_0,S)$-invariant function). Since $S^m=I$, we have
  $x(t+T)=x(t), \forall t\in \mathbb{R}$. A periodic solution $(T,x)$ of
  \eqref{s-brake-orbit} is called a $S$-{\it symmetrical brake orbit} of
  $H$. In the following we always assume $S=e^{2\pi J/m}=I_{2n}\cos \frac{2\pi}{m}+J\sin\frac{2\pi}{m}$.
\subsection{($N_0,S$)-Symmetrical Symplectic Capacity}

 We denote by $\mathcal{H}^S(U
,N_0)$ the set of $C^2$ smooth functions $H$ on $U$ satisfying the
following properties. \\
(HS1)  There is a compact set $K\subset U$ (depending on $H$) such
that $K\subset U\backslash\partial U$ and
$$H(U\backslash K)\equiv m(H) \;{\rm a \;constant}.$$
(HS2)  There is an open set $O\subset U$ and $\{0\}\in O$ (depending
on $H$) on which
$$H(O)\equiv 0.$$
(HS3)  $0\leq H(x)\leq m(H)$ for all $x\in U$.\\
(HS4)  $H(Sx)=H(N_0x)=H(x)$.

 A function $H\in \mathcal{H}^S(U,N_0)$ is called
admissible if \eqref{s-brake-orbit} has no $S$-symmetrical brake orbit, or all the $S$-symmetrical brake orbits of
\eqref{s-brake-orbit} are either constant, ie., $x(t)\equiv x(0)$,
 $\forall t\in\mathbb{R}$ or have the minimal period $T>1$. Denote the set of
admissible functions by
$\mathcal{H}^S_a(U,N_0)\subseteq\mathcal{H}^S(U,N_0)$.

\begin{defi}\label{ss}
For any  $(N_0,S)$-invariant open subset $U\subset\mathbb{R}^{2n}$, satisfying $\{0\}\in U$, the
$(N_0,S)$-symmetrical symplectic capacity is define by
\begin{equation}
c_{N_0,S}(U)=sup\{m(H)|H\in \mathcal{H}^S_a(U,N_0)\}.
\end{equation}
\end{defi}

Similar to Theorem \ref{main-result}, we have the following result.
Its proof is almost the same  as that of Theorem \ref{main-result}.
We omit the details here.
\begin{thm}\label{S-B-O} $c_{N_0,S}$ has the following properties:\\
(1) $c_{N_0,S}(U_1)\leq c_{N_0,S}(U_2)$ provided there exists
$(N_0,S)$-equivariant symplectic embedding $~~~~~\phi: (U_1,
\omega_0)\rightarrow (U_2, \omega_0)$ satisfying $\phi(N_0
x)=N_0\phi(x)$, $\phi(Sx)=S\phi(x)$ for all $x\in U_1$.\\
(2) $c_{N_0,S}(\alpha U)=\alpha^2 c_{N_0,S}(U),\alpha\neq 0,
\alpha\in
\mathbb{R}$.\\
(3) $c_{N_0,S}(B(1))=\pi,$ where $B(1)$ is the unit ball in
$\mathbb{R}^{2n}$.
\end{thm}

 \subsection{Applications for S-Symmetrical Brake Orbits}
 From Definition \ref{ss} and   Theorem \ref{S-B-O}, we can give a proof of Theorem
 \ref{T1.3} similarly as the proof of Theorem \ref{c-H-Z}. Further
 more, we have the following result.

 \begin{thm}\label{T3.2}
 If $S=e^{2\pi J/m}$, with $m\in \mathbb{N}\backslash\{1\}$, $H\in C^2(\mathbb{R}^{2n},\mathbb{R})$,
 satisfying $H(N_0z)=H(Sz)=H(z),\;\forall z\in \mathbb{R}^{2n}$  and
$H(x)\rightarrow +\infty$ with $|x|\rightarrow+\infty$, then there
exists a dense subset $O\subset(H(0),+\infty)$, such that for every
$\lambda\in O$ there is a nontrivial  $S$-symmetrical brake orbit
$(T,x(t))$ of $H$, with $H(x(t))=\lambda$.
\end{thm}
\noindent {\it Proof.} Since $H(x)\rightarrow +\infty$ with
$|x|\rightarrow+\infty$,  $H$ is bounded from below, so we can
assume $H\geq 0$. For any $M>H(0)\geq 0$ and $\varepsilon>0$, there
exist  $0<R_1<R_2$ such that
$$H(x)\leq M,\; \forall\; |x|\leq R_1,\;{\rm and}$$
$$H(x)\geq M+\varepsilon, \;\forall\; |x|\geq R_2.$$
Define a smooth function $f:\mathbb{R}^+\rightarrow \mathbb{R}$ by
$$\left\{\begin{array}{lll}
f(s)&=0,\;\;&{\rm for}\;\; s\leq M,\\
f(s)&\geq 0,\;\;&{\rm for }\;\;M< s< M+\varepsilon ,\\
f(s)&=\pi R_2^2+1,\;\;&{\rm for }\;\;s\geq M+\varepsilon,
\end{array}\right.$$
and $f'\geq 0$. Let $F(z)=f(H(z))$, then there hold
$$\left\{\begin{array}{lll}
F(z)&=0,\;\;&{\rm for}\;\; |z|\leq R_1,\\
F(z)&\geq 0,\;\;&{\rm for }\;\;R_1< |z|< R_2 ,\\
F(z)&=\pi R_2^2+1,\;\;&{\rm for }\;\;|z|\geq R_2,
\end{array}\right.$$
so $F\in\mathcal{H}^S(B(R_2+\epsilon),N_0)$ for some small
$\epsilon>0$, and  $m(F)=\pi
R_2^2+1>c_{N_0,S}(B(R_2+\varepsilon))=\pi(R_2+\varepsilon)^2$. From
Theorem \ref{S-B-O} and the Definition \ref{ss}, the following
problem
\begin{equation}
 \left\{ \begin {array}{lll}
\dot{x}(t)=J\nabla F(x(t)),\\
 x(-t)=N_0x(t),\\
x(t+\frac{T}{m})=Sx(t)
  \end{array}\right.
 \end{equation}
 has a  $T$-periodic
solution $x(t)$ with $0<T\le 1$.
 From the definition of $f$ and $F$, we have
 $$f'(H(x(t)))=\lambda > 0, \; \forall t\in \mathbb{R}.$$
 Define $y(t)=x(t/\lambda)$
 we have $(\lambda T,y(t))$ is a $S$-symmetrical brake orbit of
  $H$, and $M\leq H(y(t))\leq M+\varepsilon$.  Since $M$ and
  $\varepsilon$ are arbitrary the theorem is proved.$\hfill\Box$

\noindent{\bf Remark.} Similar to Theorem 1.2, we can prove that
every compact $(N_0,S)$-contact type hypersurface $\Sigma$ in
$\mathbb{R}^{2n}$ with $\Sigma\cap L_0\ne \emptyset$ possesses  an
$S$-symmetric closed brake characteristic. We note that the
``figure-eight orbit'' is a special case.

\end{document}